\documentclass[runningheads]{llncs}
\usepackage[T1]{fontenc}
\usepackage{graphicx}
\usepackage{xcolor}
\usepackage{amsmath}
\usepackage{amsfonts}

\begin{document}
\title{Wasserstein-enabled characterization of designs and myopic decisions in Bayesian Optimization\\
\small{\color{red}[Submitted to Special session 8: Bayesian optimization]\color{black}}}

%
%
\author{Antonio Candelieri\inst{1}\orcidID{0000-0003-1431-576X} \and
Francesco Archetti\inst{1}\orcidID{0000-0003-1131-3830}}
\authorrunning{A. Candelieri et al.}
%
\institute{ University of Milano-Bicocca, Milan, 20126, Italy
\email{\{antonio.candelieri,francesco.archetti\}@unimib.it}}
\maketitle              

\begin{abstract}
Impractical assumptions, an inherently myopic nature, and the crucial role of the initial design, all together contribute to making theoretical convergence proofs of little value in real-life Bayesian Optimization applications. In this paper, we propose a novel characterization of the design depending on its distributional properties, separately measured with respect to the coverage of the search space and the concentration around the best observed function value. These measures are based on the Wasserstein distance and enable a model-free evaluation of the information value of the design before deciding the next query. Then, embracing the myopic nature of Bayesian Optimization, we take an empirical approach to analyze the relation between the proposed characterization of the design and the quality of the next query. Ultimately, we provide important and useful insights that might inspire the definition of a new generation of acquisition functions in Bayesian Optimization.
\end{abstract}
\section{Introduction}
Bayesian Optimization (BO) \cite{archetti2019bayesian,garnett2023bayesian} is a sample-efficient sequential model-based algorithm for global optimization of black-box, expensive-to-evaluate, non-convex objective functions. The reference problem is 
\begin{equation}
    x^* \in \underset{x \in \mathcal{X}}{\arg \min} f(x)
\end{equation}
where $\mathcal{X}\subset\mathbb{R}^d$ is the so-called search space, typically box-bounded.\\

BO iterates by (\textit{i}) fitting a probabilistic surrogate model of $f(x)$ to an available set of observations $\mathcal{D}=\{(x_i,y_i)\}_{i=1:n}$, usually called \textit{design}, and (\textit{ii}) selecting the next \textit{query}, namely $x'$, by optimizing an acquisition function that balances between exploration and exploitation, according to the surrogate.

BO is known to be \textit{sample efficient}, meaning that it finds a better solution than other global optimization strategies, within a few queries. Thanks to this property, it has been and is still widely adopted to address industrial-relevant applications \cite{li2026bayesian,siska2026guide} and also Machine Learning tasks \cite{alcobacca2026literature,elshewey2026enhancing,hiccyilmaz2026bayesian}.

Although theoretical convergence proofs exist for specific BO implementations, such as the widely known \cite{srinivas2012information}, actual convergence in real-life settings is a critical issue, mostly due to impractical assumptions of convergence theorems. 

The most relevant and well-known obstacle to global convergence of BO is the mispecification of the probabilistic surrogate model, typically a Gaussian Process (GP) \cite{gramacy2020surrogates,williams2006gaussian}. Theoretical convergence proofs hold only in the standard setting \cite{bogunovic2021misspecified}, meaning that the objective function is a member of a specific Reproducing Kernel Hilbert Space (RKHS) with bounded norm, known a-priori. In other terms, convergence to the optimum is proved under the assumption that the GP kernel's hyperparameters are known in advance, which is not the case in practice. In BO, the GP kernel's hyperparameters are tuned depending on the available design $\mathcal{D}$ via Maximum Likelihood Estimation (MLE) or Maximum-A-Posteriori (MAP), leading common "no regret" acquisition functions to get stuck in local optima \cite{bull2011convergence,candelieri2024mle}. In simpler terms, sequential queries are model-based and thus hyperparameter-sensitive; the acquisition function can introduce a feedback loop where sequential queries reinforce bad hyperparameters.
Heuristic methods \cite{wang2014theoretical,wang2016bayesian,wabersich2016advancing,berkenkamp2019no} have been proposed over time, basically using Lemma 4 of \cite{bull2011convergence}, stating that decreasing the \textit{length-scale} hyperparameter of the kernel leads to consider a larger RKHS, increasing exploration at the expense of a higher regret. Other methods to increase the chance of global convergence propose random exploration steps  \cite{candelieri2023mastering,kim2025enhancing} or inexact optimization of the acquisition function \cite{kim2025bayesian}.

Another issue, strictly related to model mispecification, is that the random design generated to initialize BO can definitely prevent from reaching the global optimum, regardless of the acquisition function and its specific exploitation-exploration trade-off mechanism. Indeed, the typical experimental setting for BO consists in averaging on multiple independent runs from different initial designs, all with the same number of random observations, aiming at mitigating the impact of random initialization. However, this way of proceeding is impractical for many real-life problems with highly expensive objective functions.

Indeed, many recent papers propose strategies to dynamically select the behavior of the acquisition function \cite{candelieri2023mastering,kim2025enhancing}, the model's hyperparameters \cite{wang2014theoretical,wang2016bayesian,wabersich2016advancing,berkenkamp2019no}, or both \cite{park2025boost}, to possibly revert from being stack into local optima. The long-term perspective proposed in \cite{bossek2020initial}, 
consists of developing approaches that dynamically decide whether to take the next sample from a (quasi-)random distribution or whether to derive it from the surrogate model, with an expected substantial improvement of performance, especially when the maximum allowed number of queries is known in advance.

A specific analysis on the difficulty of generating a good design for GP fitting -- not necessarily for BO -- is given in \cite{zhang2021distance}. While \textit{maximin} designs can be awful, Latin Hypercube Sampling (LHS) methods offer better performance. However, neither method is well suited for GP's hyperparameters inference much so that random design can surprisingly result in the most effective strategy. Specifically for (few-shot) BO, \cite{hvarfner2025informed} has recently proposed a novel acquisition function that balances predictive uncertainty reduction with hyperparameter learning using information-theoretic principles.

Finally, another -- often overlooked -- issue affecting the practical convergence of BO is that (almost all) the acquisition functions are \textit{myopic}. Knowledge Gradient (KG) \cite{frazier2008knowledge} and Predictive Entropy Search (PES) \cite{hernandez2014predictive} are exceptions, which have been shown to be effective but entailing a relevant computational overhead for nested integrations or Monte Carlo sampling.\\

Motivated by the issues discussed so far, this paper aims to contribute to this challenge as follows:
\begin{itemize}
    \item proposing a novel and statistically principled mechanism, based on the Wasserstein distance, to characterize -- aka \textit{featurize} -- a design depending on its distributional properties; 
    \item analyzing, according to the proposed characterization, the differences between designs suitably generated for representing typical situation in BO (i.e., LHS, queries converging to the global optimum, and queries converging to different a solution);   
    \item investigating possible relations between the characterization of a design, the risk of model mispecification, and the possible improvement provided by different myopic acquisition functions.
\end{itemize}\color{black}
For the empirical analysis, a set of well-diversified test problems has been considered. The related works are the references quoted so far.

\section{Background}

\subsection{Bayesian Optimization in brief}
Given a set of available observations $\mathcal{D}=\{(x_i,y_i)\}$, we separately denote the set of queried locations by $X=\{x_i\}_{i=1:n}$ and the set of associated values by $Y=\{y_i=f(x_i)\}_{i=1:n}$, with $f(x)$ assumed \textit{noise-free} in this paper.

As a probabilistic surrogate model, we consider a GP, with prediction and associated (model-based) uncertainty denoted by $\mu(x|\mathcal{D})$ and $\sigma(x|\mathcal{D})$, respectively. For the sake of readability, we omit their well-known equations.

As far as the acquisition function is concerned, here, we introduce a slightly different notation from usual, and consider that the next location to query $x'$ is given by a parametric \textit{policy} $\pi(\mathcal{D|\theta})$, consisting of GP fitting followed by the optimization of a specific acquisition function. Formally,
\begin{equation}
    x' \leftarrow \pi(\mathcal{D}|\theta)
\end{equation}
where $\theta$ denotes all the hyperparameters involved, both the GP's and the acquisition function's ones.

Embracing the myopic nature of common acquisition functions, the goal of $\pi(\mathcal{D|\theta})$ should be to obtain $x'=x^*$ or, more generally, to provide $x'$ such that $x^* \in \mathcal{D}\cup \{(x',y')\}$, with $y'=f(x')$, in order to consider both the case in which the optimum is already in $\mathcal{D}$ or it is exactly $x'$.

\subsection{Wasserstein distance in brief}
The Wasserstein distance comes from the Optimal Transport (OT) theory \cite{peyre2019computational,santambrogio2015optimal} and -- contrary to divergence measures -- it is a \textit{distance metric} over the space of probability distributions. Moreover, it can be used to compare two discrete or two continuous probability distributions, as well as one discrete and one continuous.
Due to these properties, it has recently gained a central role in Machine Learning \cite{montesuma2024recent} and Generative AI \cite{cheng2024convergence,gao2025wasserstein}.
For a detailed overview on OT and the Wasserstein distance, the reader can refer to \cite{peyre2019computational,santambrogio2015optimal}; here, we summarize the background needed for the scope of this paper.\\

Our setting consists only of empirical discrete probability distributions. Let $\alpha$ and $\beta$ be two of that kind, we can denote them by:
\begin{equation}
    \alpha=\frac{1}{m_\alpha}\sum_{j=1}^{m_\alpha}\delta_{a_j}\;\;, \;\; \beta=\frac{1}{m_\beta}\sum_{l=1}^{m_\beta}\delta_{b_l}
\end{equation}
with $a_j,b_l\in\mathbb{R}^d$ and $\delta_w$ denoting the Dirac's delta function centered at $w$.

Another way to look at two empirical discrete probability distributions is to consider them as point clouds, represented as matrices: $A\in\mathbb{R}^{m_\alpha\times d}$, with $A_{j*}=a_j$, and $B\in\mathbb{R}^{m_\beta\times d}$, with $B_{l*}=b_j$. In other terms, $A$ and $B$ are the discrete supports of the two uniform distributions $\alpha$ and $\beta$, respectively.

Computing the 2-Wasserstein distance between $\alpha$ and $\beta$ means solving the following linear programming problem (Kantorovich's formulation):
\begin{equation}
\begin{aligned}
    \mathcal{W}_2(\alpha,\beta) = \underset{T \in \mathbb{R}^{m_\alpha\times m_\beta}}{\min}\; & \left[\sum_{\substack{j=1:m_\alpha \\ l=1:m_\beta}} \|a_j - b_l\|^2 T_{jl}\right]^\frac{1}{2}\\
    s.t.\; & T^\top\mathbf{1}_{m_\alpha} = B\\
    & T\mathbf{1}_{m_\beta} = A\\
    & T_{ij} \geq 0,\;\forall i=1:m_\alpha, j=1:m_\beta
\end{aligned}
\end{equation}

\noindent
where $\mathbf{1}_{q}$ denotes a vector of all-ones of $q$ dimensions.\\

In the rest of the paper, we use $\mathcal{W}_2(A,B)$ or $\mathcal{W}_2(\alpha,\beta)$ interchangeably, depending on which notation is more convenient for the discussion.

\section{Information value of $\mathcal{D}$ \textit{vs} quality of myopic decisions}

Our research hypothesis is that the information value of the current design $\mathcal{D}$ can be described by the distributional properties of its two sets $X$ and $Y$, which can be jointly analyzed to predict the quality of the next query $x'$, \textbf{without relying on any specific model}. 

First, we assume that there exists a vector-valued function $\mathcal{S}:[\mathcal{X}\times \mathbb{R}]^n\rightarrow \mathbb{R}^h$ characterizing the information value of $\mathcal{D}=(X,Y)$ with respect to $h$ different features. In this study, we propose $h=2$, with $\mathcal{S}_1(\mathcal{D})$ and $\mathcal{S}_2(\mathcal{D})$ denoting the information value of $\mathcal{D}$ with respect to the coverage of the search space $\mathcal{X}$ and the observed variability of $f(x)$, respectively.

\subsection{$\mathcal{S}_1(\mathcal{D})$: information value of $\mathcal{D}$ with respect to $\mathcal{X}$}
$\mathcal{S}_1(\mathcal{D})$ is computed as the 2-Wasserstein distance between the available set of queried locations, namely $X$, and a regular $m\times m$ grid $G_{\mathcal{X}}=\{g_j\}_{j=1:m^2}$ over the search space $\mathcal{X}$. Formally,
\begin{equation}
    \mathcal{S}_1(\mathcal{D}) \stackrel{\text{def}}{=} \mathcal{W}_2^2(X,G_\mathcal{X})
\end{equation}

Thus, $\mathcal{S}_1(\mathcal{D})$ quantifies how far the set $X$ is from the uniform distribution approximated by $G_\mathcal{X}$. A value of $\mathcal{S}_1(\mathcal{D})$ close to zero means that $X$ is spread almost uniformly over $\mathcal{X}$.

Obviously, any discrepancy measure can be alternatively used to quantify $\mathcal{S}_1(\mathcal{D})$, but we suggest to choose among those which are distance metrics (e.g., Star Discrepancy, Maximum-Mean-Discrepancy, Total Variation) and to avoid divergences, such as the Kullback-Leibler divergence.

We use the 2-Wasserstein distance to ensure coherence throughout the proposed framework, since it is the only suitable choice to measure the information value of $Y$, as explained in the next section.

\subsection{$\mathcal{S}_2(\mathcal{D})$: information value of $\mathcal{D}$ with respect to $f(x)$}
Since the objective function $f(x)$ is black-box, its co-domain is unknown, and, contrary to $\mathcal{S}_1(\mathcal{D})$, it is not possible to pre-define a regular grid of values to be compared against the available set $Y$.

On the other hand, an obvious property of a good $Y$ is to show some variability, since the trivial assumption is that $f(x)$ is not flat. Although the standard deviation of $Y$ is an intuitive choice, it is not suitable to quantify the difference between various distributions of $Y$ and their concentration near the best value observed so far, namely the \textit{best seen} $y^+=\min\{Y\}$. Thus, we use the 2-Wasserstein distance between $Y$ and a Dirac's delta over $y^+$:
\begin{equation}
    \mathcal{S}_2(\mathcal{D}) \stackrel{\text{def}}{=} \mathcal{W}_2^2(Y,\delta_{y^+})
\end{equation}

\noindent
It is trivial to prove that $\mathcal{W}_2^2(Y,\delta_{y^+})=\sum_{i=1}^n (y_i-y^+)=\big(\sum_{i=1}^n y_i\big)-ny^+$. \\

Finally, any $\mathcal{D}$ can be represented as a point in the space spanned by $\mathcal{S}_1(\mathcal{D})$ and $\mathcal{S}_2(\mathcal{D})$, which are the distributional properties of the sets $X$ and $Y$.

\begin{figure}
    \centering
    \includegraphics[width=0.65\textwidth]{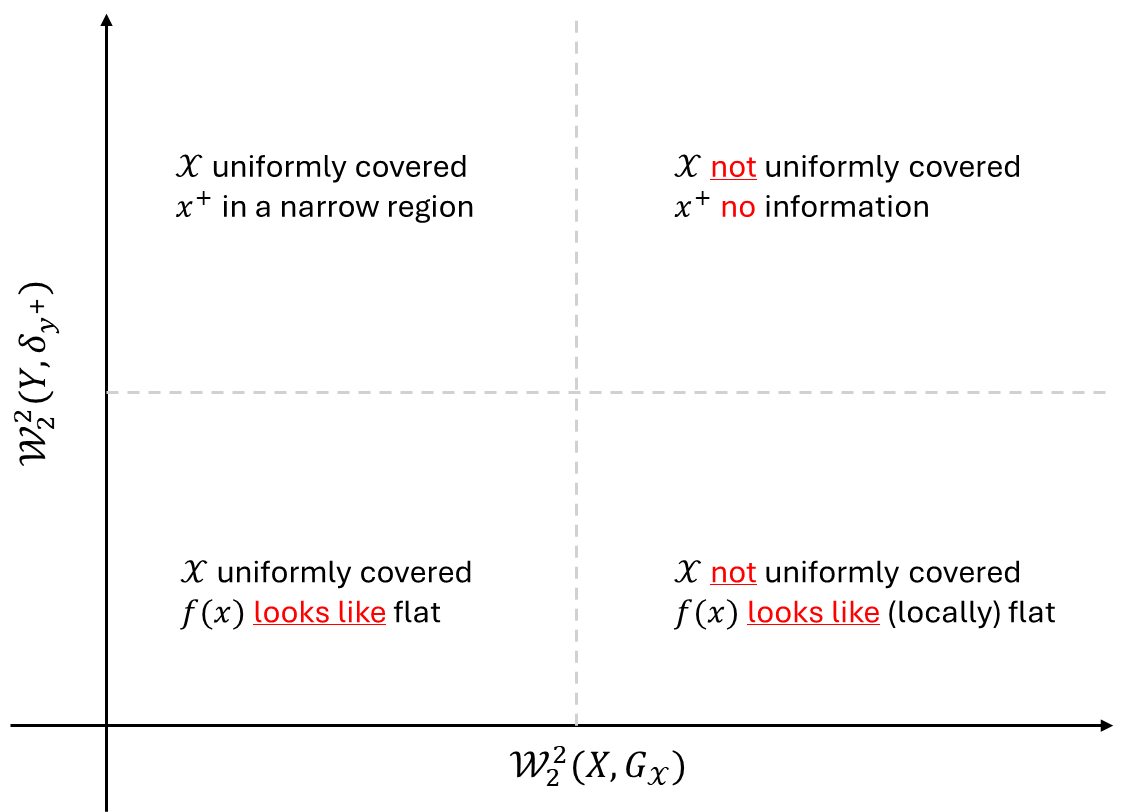}
    \caption{A conceptual representation of how a design $\mathcal{D}$ is characterized by its distributional properties $\mathcal{S}_1(\mathcal{D})\stackrel{\text{def}}{=}\mathcal{W}_2^2(X,G_\mathcal{X})$ and $\mathcal{S}_2(\mathcal{D})\stackrel{\text{def}}{=}\mathcal{W}_2^2(Y,\delta_{y^+})$.}
    \label{fig:schema}
\end{figure}

Figure \ref{fig:schema} provides a graphical representation of how the location of a given $\mathcal{D}$ within this space can provide useful insights about the information value of $\mathcal{D}$ itself. 
A small $\mathcal{W}_2^2(X,G_\mathcal{X})$ indicates that $X$ covers the search space $\mathcal{X}$ almost uniformly. A small $\mathcal{W}_2^2(Y,\delta_{y^+})$ suggests that $f(x)$ "looks like" flat: we remark that it "looks like" flat because $Y$ depends on the locations in $X$. More important, a large $\mathcal{W}_2^2(Y,\delta_{y^+})$ and a small $\mathcal{W}_2^2(X,G_\mathcal{X})$ suggest that the best solution of the design, denoted by $(x^+,y^+)$, could lie within a narrow region.
\vspace{-0.5cm}

\subsection{Model mispecification and quality of myopic decisions}
While the previously defined $\mathcal{S}_1(\mathcal{D})$ and $\mathcal{S}_2(\mathcal{D})$ are fully specified before the myopic decision $x'$, here we define two measures aimed at evaluating a possible model mispecificaton implied by $\mathcal{D}$ -- under the adopted model fitting procedure (e.g., MLE) -- and the quality of the resulting myopic decision.

The first measure is simply the root mean squared error (RMSE) between $f(x)$ and the GP's prediction for the points of the grid $G_\mathcal{X}$, that is:
\begin{equation}
    RMSE = \sqrt{{\frac{1}{\vert G_\mathcal{X}\vert}\sum_{x_i\in G_\mathcal{X}} \Big(f(x_i)-\mu(x_i|\mathcal{D},\theta) \Big)^2 }}
\end{equation}

Clearly, this measure cannot be computed in a real-life setting; it is used in this study just to quantify GP's mispecification in our experiments.

The second measure, denoted by $\Delta_y$, is the improvement/worsening implied by $y'=f(x')$ with respect to $y^+$, formally:
\begin{equation}
    \Delta_y = y^+-y'
\end{equation}
where a value of $\Delta_y$ lower/higher than zero indicates that $y'$ is worse/better than $y^+$, while $\Delta_y$ is equal to zero when $y'=y^+$. Unlike RMSE, $\Delta_y$ can be computed in a real-life BO task, but only after $x'$ is evaluated. Furthermore, also the immediate regret $y'-y^*$ is  computed but, similarly to RMSE, only for analytical purposes.

\section{Experiments and results}

\subsection{Setting}
In this section, we summarize the relevant details of our experimental setting:
\begin{itemize}
    \item eight test problems with $d=1$ and six with $d=2$ (all search spaces have been preliminary rescaled to $[0,1]^d$). A graphical representation of the test problems is provided in Appendix;
    \item three types of design $\mathcal{D}$, generated as follows:
        \begin{itemize}
            \item (1) from LHS;
            \item (2) 50\% from LHS and 50\% from a neighborhood $N(x^*)$ of the true optimizer to simulate a BO process converging to $x^*$;
            \item (3) 50\% from LHS and 50\% from the neighborhood $N(\tilde{x})$ of a random point $\tilde{x}\neq x^*$ to simulate a BO process converging to a wrong solution;
        \end{itemize}    
    \item number of observations in $\mathcal{D}$: $n \in \{5d,\lfloor12.5d\rfloor,20d\}$;
    \item GP model: Matérn$_{3/2}$ kernel; fitting via MLE and nugget estimation only in the case of an ill-conditioned kernel matrix; no standardization/scaling of $Y$;
    \item Acquisition function: to avoid randomness implied by the method used to optimize the acquisition function, we just pick the best solution on the finite grid $G_\mathcal{X}$ (consisting of 10.000 points). Four acquisition functions, all sharing the same fitted GP, are considered:
    \begin{itemize}
        \item surface-response (SR): $\pi(x|\mathcal{D},\theta) = \underset{x \in G_\mathcal{X}}{\arg \min}\; \mu(x|\mathcal{D})$
        \item uncertainty reduction -- i.e., maximization of the GP's standard deviation (SD): $\pi(x|\mathcal{D},\theta) = \underset{x \in G_\mathcal{X}}{\arg \max}\; \sigma(x|\mathcal{D})$ 
        \item Expected Improvement (EI): $\pi(x|\mathcal{D},\theta)= \underset{x \in G_\mathcal{X}}{\arg \max}\; EI(x|\mathcal{D})$
        \item Lower Confidence Bound (LCB): $\pi(x|\mathcal{D},\theta)= \underset{x \in G_\mathcal{X}}{\arg \min}\; \mu(x|\mathcal{D}) - \beta\sigma(x|\mathcal{D})$, with $\beta=1$;
    \end{itemize}    
\end{itemize}

To guarantee reproducibility of the experiments, the code is available at the following Github repository:\\
\url{https://github.com/acandelieri/LION20\_Candelieri\_Archetti.git}

\subsection{Results on characterizing designs through $\mathcal{S}_1(\mathcal{D})$ and  $\mathcal{S}_2(\mathcal{D})$ }
In this section, we summarize the most relevant results on the characterization of different designs (i.e., type and size) on a set of diversified test problems. The characterization is given by the two proposed functions $\mathcal{S}_1(\mathcal{D})$ and $\mathcal{S}_2(\mathcal{D})$.

\begin{figure}[h!]
    \centering
    \includegraphics[width=0.7\textwidth]{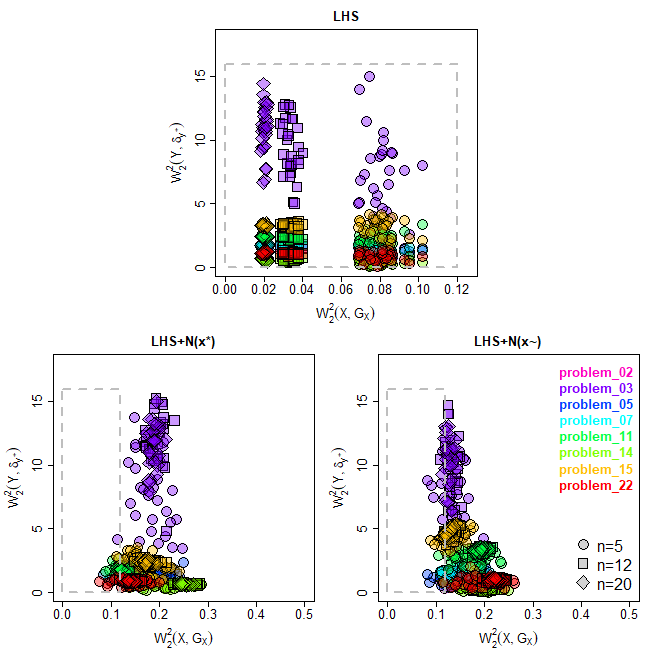}
    \caption{Thirty designs, with size 5, 12, and 20, on eight 1-dimensional test problems, separately for three types of design: LHS, LHS+$N(x^*)$, and LHS+$N(\tilde{x})$.}
    \label{fig:res1}
\end{figure}

\begin{figure}[h!]
    \centering
    \includegraphics[width=0.7\textwidth]{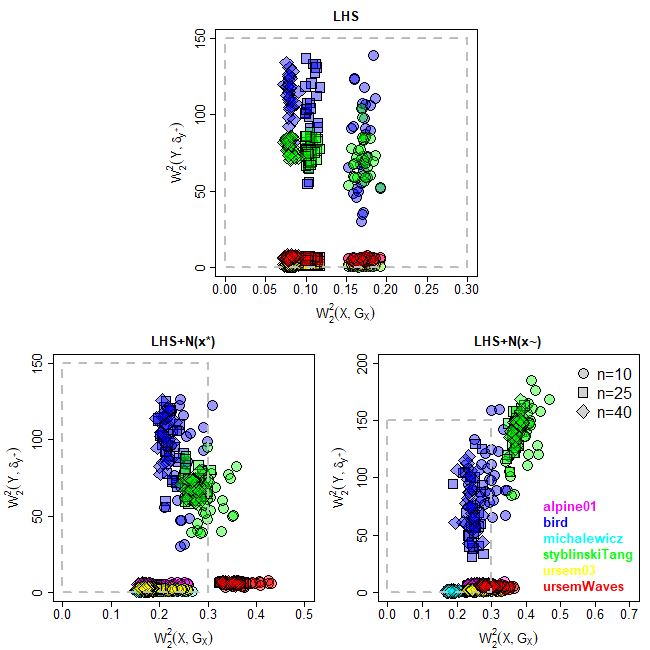}
    \caption{Thirty designs, with size 10, 25, and 40, on six 2-dimensional test problems, separately for three types of design: LHS, LHS+$N(x^*)$, and LHS+$N(\tilde{x})$.}
    \label{fig:res2}
\end{figure}

Figure \ref{fig:res1} and Figure \ref{fig:res2} refer to the 1-dimensional and the 2-dimensional test problems, respectively. Each figure is organized into three charts, one for each type of design, while different shapes and colors of the points refer to the size of $\mathcal{D}$ and the specific test problem, respectively. Since LHS leads to significantly lower values of $\mathcal{W}_2^2(X.G_\mathcal{X})$ than the other two design types, the axis ranges of the three charts are different; the dashed gray rectangle helps for the comparison.

\newpage

\subsubsection{As far $\mathcal{S}_1(\mathcal{D})\stackrel{\text{def}}{=}\mathcal{W}_2^2(X,G_\mathcal{X})$ is concerned,} we observe that it clearly depends on $n$: the higher the value of $n$, the lower the 2-Wasserstein distance between $X$ and $G_\mathcal{X}$. Furthermore, the deviation in the values of $\mathcal{W}_2^2(X,G_\mathcal{X})$ also decreases as $n$ increases. Finally, the values of $\mathcal{W}_2^2(X,G_\mathcal{X})$ are significantly lower than those obtained with $\mathcal{D}$ from LHS+$N(x^*)$ and LHS+$N(\tilde{x})$.

This "well-shaped" organization of the designs in the $\mathcal{S}_1(\mathcal{D})\times\mathcal{S}_2(\mathcal{D})$ space -- with respect to the x-axis -- has been observed in both the 1-dimensional and 2-dimensional test problems, and is in line with the expectation about the LHS technique: the design becomes more uniformly distributed in $\mathcal{X}$ with $n$ increasing.

\subsubsection{As far as $\mathcal{S}_2(\mathcal{D})\stackrel{\text{def}}{=}\mathcal{W}_2^2(Y,\delta_{y^+})$ is concerned,} there is no clear relation with $n$. Indeed, $\mathcal{S}_2(\mathcal{D})$ depends on $f(x)$ together with $\mathcal{D}$ and is uncorrelated with $\mathcal{W}_2^2(X,G_\mathcal{X})$, regardless of the type of design. On the other hand, an important insight is that when $f(x)$ is highly "oscillating", with many local optima, the associated $\mathcal{W}_2^2(Y,\delta_{y^+})$ is quite high.

\newpage

\subsection{Results on the quality of $x'$}
\begin{figure}[h!]
    \centering
    \includegraphics[width=\textwidth]{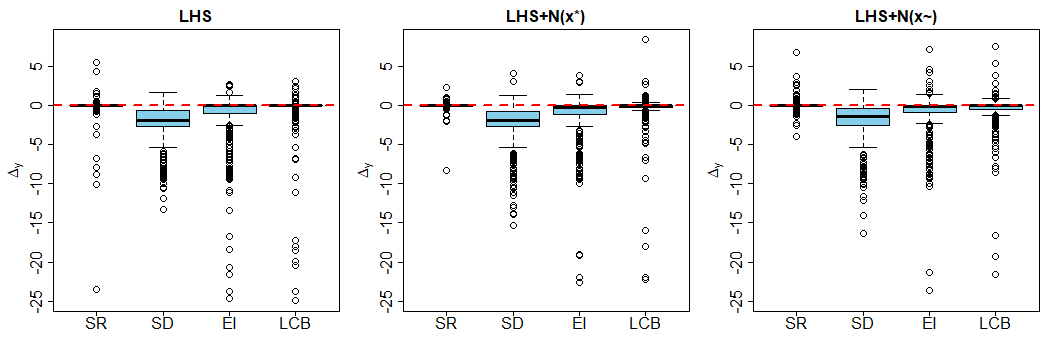}
    \caption{Boxplot of the $\Delta_y$ values on the \textbf{1-dimensional problems}: one chart for each type of design. Boxes refer to acquisition functions (i.e., SR: surface response, SD: maximization of the GP's standard deviation, EI: Expected Improvement, and LCB: Lower Confidence Bound). All the acquisition functions share the same GP model.}
    \label{fig:res3}
\end{figure}

This section summarizes the results related to the quality of the decision $x'$. Figure \ref{fig:res3} shows the differences in terms of $\Delta_y$ for the 1-dimensional test problems: a boxplot for each type of design and a box for each of the four acquisition functions. It is important to recall that the acquisition functions share the same GP, so that a possible mispecification would affect all of them. All 1-dimensional problems and all sizes $n=\{5,12,20\}$ of the design are considered. Surprisingly, the simple SR is the most robust acquisition function for all three types of design, with LCB and EI providing a more negative $\Delta_y$ than SD in some cases.

\begin{figure}[h!]
    \centering
    \includegraphics[width=\textwidth]{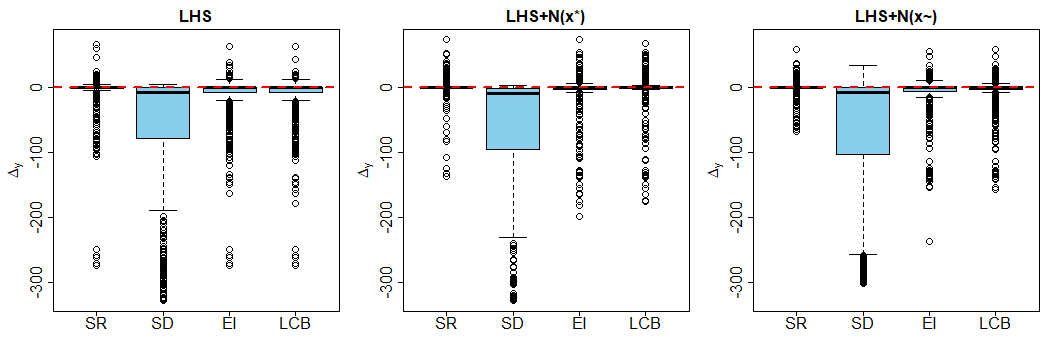}
    \caption{Boxplot of the $\Delta_y$ values on the \textbf{2-dimensional problems}: one chart for each type of design. Boxes refer to acquisition functions (i.e., SR: surface response, SD: maximization of the GP's standard deviation, EI: Expected Improvement, and LCB: Lower Confidence Bound). All the acquisition functions share the same GP model.}
    \label{fig:res4}
\end{figure}

Figure \ref{fig:res4} reports the same kind of analysis, but for the 2-dimensional test problems. The previous considerations are confirmed: SR is the most robust acquisition function, SD does not lead to any improvement, and bot LCB and EI provide more negative values of $\Delta_y$ than SR.

\newpage

\begin{figure}[h!]
    \centering
    \includegraphics[width=\textwidth]{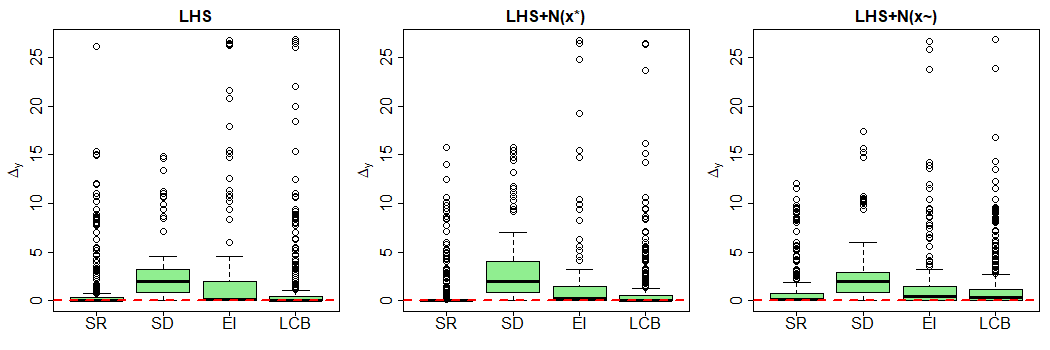}
    \caption{Boxplot of the immediate regret $y'-y^*$ on the \textbf{1-dimensional problems}.}
    \label{fig:res5}
\end{figure}

\begin{figure}[h!]
    \centering
    \includegraphics[width=\textwidth]{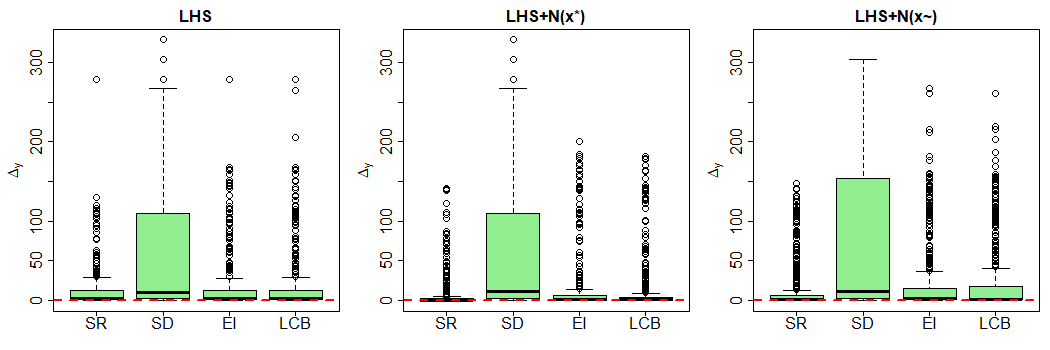}
    \caption{Boxplot of the immediate regret $y'-y^*$ on the \textbf{2-dimensional problems}.}
    \label{fig:res6}
\end{figure}

Even more interesting, SR results in the most robust choice also in terms of immediate regret $y'-y^*$, as depicted in Figure \ref{fig:res5} and Figure \ref{fig:res6}.

\newpage

It is interesting to note that this result is consistent with the empirical evidence reported in \cite{benjamins2022pi}, showing that switching to an exploration-biased acquisition function such as Probability of Improvement, close to the end of the BO process, can outperform BO based on a single no-regret acquisition function.

\subsection{Relations between $\mathcal{S}_1(\mathcal{D}),\mathcal{S}_2(\mathcal{D})$, RMSE, and $\Delta_y$}
Here, we restrict our attention to the case $n=20d$, given the empirical evidence suggesting that BO can find an optimal solution within 15-20 time the number of dimensions. From our experiments, we observed that there exists a (linear) correlation between $\mathcal{S}_2(\mathcal{D})\stackrel{\text{def}}{=}\mathcal{W}_2^2(Y,\delta_{y^+})$ and GP's mispecification as measured by RMSE. Table \ref{tab:1} reports the Pearson's correlation coefficient -- and the associated $p$-value -- separately for the type of designs (by row) and the dimensionality of the test problems (by column).

In simpler terms, the fewer/more values in $Y$ are concentrated close to $y^+$ (i.e., large/small values of $\mathcal{W}_2^2(Y,\delta_{y^+})$) the higher/lower the risk of model mispecification. This holds for the three types of design, with high correlation values, even if slightly different from a design type to another. The most important insight is that LHS cannot avoid the risk of model mispecification \cite{hvarfner2025informed,zhang2021distance}.

\newpage

\begin{table}[h!]
    \centering
    \begin{tabular}{|l|c|c|}
        \hline
        \;design type\; & \;1-dimensional test problems\; & \;2-dimensional test problems\; \\
        \hline\hline
        \;LHS\; & 0.9377 \scriptsize{(\textit{p}-value<0.0001)} &  0.8763 \scriptsize{(\textit{p}-value<0.0001)}\\
        \;LHS+$N(x^*)$\; & 0.9644 \scriptsize{(\textit{p}-value<0.0001)} & 0.9331 \scriptsize{(\textit{p}-value<0.0001)} \\
        \;LHS+$N(\tilde{x})$\; & 0.8631 \scriptsize{(\textit{p}-value<0.0001)} & 0.6891 \scriptsize{(\textit{p}-value<0.0001)} \\
        \hline
    \end{tabular}
    \caption{Pearson's correlation coefficient between $\mathcal{W}_2^2(Y,\delta_{y^+})$ and RMSE, $n=20$.}
    \label{tab:1}
\end{table}

\vspace{-1cm}

Based on the results on the quality of $x'$ for the different acquisition functions, we consider only the SR acquisition function and report the correlations between RMSE and $\Delta_y^{SR}$ (Table \ref{tab:2}). The most interesting result refers to LHS, with a (statistically significant) negative correlation -- even if weak -- which has been observed on both the 1-dimensional and 2-dimensional test problems. This holds for LHS, but not for the other two design types, especially for LHS+$N(x^*)$ because, even if the GP can be globally mispecified, it is hardly mispecified close to $x^*$ and thus SR has more chances to improve with respect to $y^+$.

\vspace{-0.5cm}

\begin{table}[h!]
    \centering
    \begin{tabular}{|l|c|c|}
        \hline
        \;design type\; & \;1-dimensional test problems\; & \;2-dimensional test problems\; \\
        \hline\hline
        \;LHS\; & -0.2526 \scriptsize{(\textit{p}-value<0.0001)} & -0.4192 \scriptsize{(\textit{p}-value<0.0001)} \\
        \;LHS+$N(x^*)$\; & 0.2521 \scriptsize{(\textit{p}-value<0.0001)} & 0.1285 \scriptsize{(\textit{p}-value=0.086)} \\
        \;LHS+$N(\tilde{x})$\; & 0.2538 \scriptsize{(\textit{p}-value<0.0001)} & -0.1616 \scriptsize{(\textit{p}-value=0.010)} \\
        \hline
    \end{tabular}
    \caption{Pearson's correlation coefficient between RMSE and $\Delta_y^{SR}$, $n=20d$.}
    \label{tab:2}
\end{table}

\vspace{-1cm} 
The most important result comes from an obvious question: since $\mathcal{W}_2^2(Y,\delta_{y^+})$ can be calculated before $x'$ and is correlated with RMSE, and since RMSE is inversely correlated with $\Delta_y^{SR}$, is there a correlation between $\mathcal{W}_2^2(Y,\delta_{y^+})$ and $\Delta_y^{SR}$? Luckily, the answer to this question is \textit{"yes, for LHS and SR"}: Table \ref{tab:3} reports the correlation coefficient between $\mathcal{W}_2^2(Y,\delta_{y^+})$ and $\Delta_y^{SR}$. In the case of an LHS design, we again observe a (statistically significant) negative correlation, this time between $\mathcal{W}_2^2(Y,\delta_{y^+})$ and $\Delta_y^{SR}$, on both 1-dimensional and 2-dimensional test problems. For the other two types of design, the correlation is positive but statistically significant only for the 1-dimensional problems.

\vspace{-0.5cm}

\begin{table}[h!]
    \centering
    \begin{tabular}{|l|c|c|}
        \hline
        \;design type\; & \;1-dimensional test problems\; & \;2-dimensional test problems\; \\
        \hline\hline
        \;LHS\; & -0.2087 \scriptsize{(\textit{p}-value<0.0001)}  & -0.2569 \scriptsize{(\textit{p}-value<0.0001)}  \\
        \;LHS+$N(x^*)$\; & 0.3009 \scriptsize{(\textit{p}-value<0.0001)}  & 0.1280 \scriptsize{(\textit{p}-value=0.087)}\\
        \;LHS+$N(\tilde{x})$\; & 0.2295 \scriptsize{(\textit{p}-value<0.0001)}  & 0.0556 \scriptsize{(\textit{p}-value=0.458)}\\
        \hline
    \end{tabular}
    \caption{Pearson's correlation coefficient between $\mathcal{W}_2^2(Y,\delta_{y^+})$ and $\Delta_y^{SR}$, $n=20d$.}
    \label{tab:3}
\end{table}

\vspace{-1cm}
Finally, we have always observed a statistically significant linear positive correlation between $\mathcal{W}_2^2(Y,\delta_{y^+})$ and immediate regret, with a smaller correlation coefficient for SR and LCB and a larger one for EI and SD. This further remarks that high values of $\mathcal{W}_2^2(Y,\delta_{y^+})$ indicate a more complicated to be optimized $f(x)$.

\section{Discussions on results and opportunities}
We summarize the most important insights from our empirical study:
\begin{enumerate}
    \item the risk of model mispecification increases with $\mathcal{S}_2(\mathcal{D})\stackrel{\text{def}}{=}\mathcal{W}_2^2(Y,\delta_{y^+})$ and regardless of $\mathcal{S}_1(\mathcal{D})\stackrel{\text{def}}{=}\mathcal{W}_2^2(X,G_\mathcal{X})$. Thus, even if the design $\mathcal{D}$ in a generic BO iteration is close to LHS, this does not prevent model mispecification.
    \item if $\mathcal{D}$ is close to LHS and $\mathcal{S}_2(\mathcal{D})\stackrel{\text{def}}{=}\mathcal{W}_2^2(Y,\delta_{y^+})$ is "small", then a simple SR acquisition function can provide a high $\Delta_y^{SR}$, which is hopefully positive. The main issue -- at the moment -- is to quantify what "small" means.
    \item if $\mathcal{D}$ is a design associated with a BO task that is converging towards a solution -- be it $x^*$ or not -- and $\mathcal{S}_2(\mathcal{D})\stackrel{\text{def}}{=}\mathcal{W}_2^2(Y,\delta_{y^+})$ is "large", again a simple SR acquisition function can offer the chance of high $\Delta_y^{SR}$.
    \item the remaining two cases are the most critical: (a) $\mathcal{D}$ is close to LHS and $\mathcal{S}_2(\mathcal{D})\stackrel{\text{def}}{=}\mathcal{W}_2^2(Y,\delta_{y^+})$ is "large", or (b) $\mathcal{D}$ is a design converging to a solution and $\mathcal{S}_2(\mathcal{D})\stackrel{\text{def}}{=}\mathcal{W}_2^2(Y,\delta_{y^+})$ is "small". In these cases, no one of the acquisition functions considered has a relevant chance of leading to $\Delta_y\geq0$, but due to two completely different reasons: $f(x)$ is "complicated" (a) and the BO is (or is going to be) stuck in a local/global optimum (b).
\end{enumerate}

The last point of this list represents the first relevant opportunity that can follow from this paper. For the two cases described in 4., a new policy $\hat{\pi}(\mathcal{D}\vert\theta)$ should be defined with the aim of resorting $\mathcal{D}$ to one of the two more favorable situations in 2. and 3. --- i.e., "moving" $\mathcal{D}$ towards LHS if $\mathcal{S}_2(\mathcal{D})\stackrel{\text{def}}{=}\mathcal{W}_2^2(Y,\delta_{y^+})$ is "small" or forcing convergence towards a solution if it is "large".

Thus, a new generation acquisition function should switch between $\hat{\pi}(\mathcal{D\vert\theta})$ and $\pi(\mathcal{D}\vert\theta)$ (i.e., SR) jointly depending on the values of $\mathcal{S}_1(\mathcal{D})\stackrel{\text{def}}{=}\mathcal{W}_2^2(X,G_\mathcal{X})$ and $\mathcal{S}_2(\mathcal{D})\stackrel{\text{def}}{=}\mathcal{W}_2^2(Y,\delta_{y^+})$, \textbf{before any model fitting}. 

\section{Conclusions}
This paper presents a new characterization of a design based on its distributional properties measured through the Wasserstein distance. Empirical analysis of myopic decisions induced by common basic acquisition functions led to the observation of a relationship between the proposed characterization of the design, immediate regret, and improvement with respect to the best observed value in the design. This relation suggests a way for dynamically changing the behavior of the acquisition step, more principled with respect to a simple switch towards more/less exploration. Indeed, instead of just increasing exploration, queries should be selected to improve distributional properties of the design up to an optimal myopic decision can be made. It would also be interesting to investigate whether the distributional properties of the design could be related to stopping criteria, such as the Gittin's index proposed in \cite{scully2025gittins}.

Although the empirical results of this study are interesting, limitations regard a suitable quantification of what "small" and "large" $\mathcal{W}_2^2(Y,\delta_{y^+})$ really mean. Currently, this is the main challenge that our research is now focusing on in order to design a new generation of acquisition functions in BO.

\begin{credits}
\subsubsection{\ackname} This work has benefitted from the first author's participation in Dagstuhl Seminar 25451
Bayesian Optimisation (Nov 02 – Nov 07, 2025).

\subsubsection{\discintname}
The authors have no competing interests to declare that are
relevant to the content of this article.
\end{credits}

\newpage

%
%
%
\bibliographystyle{splncs04}
\bibliography{LION20_Candelieri}

\clearpage
\section{Appendix}
\begin{figure}
    \centering
    \includegraphics[width=\linewidth]{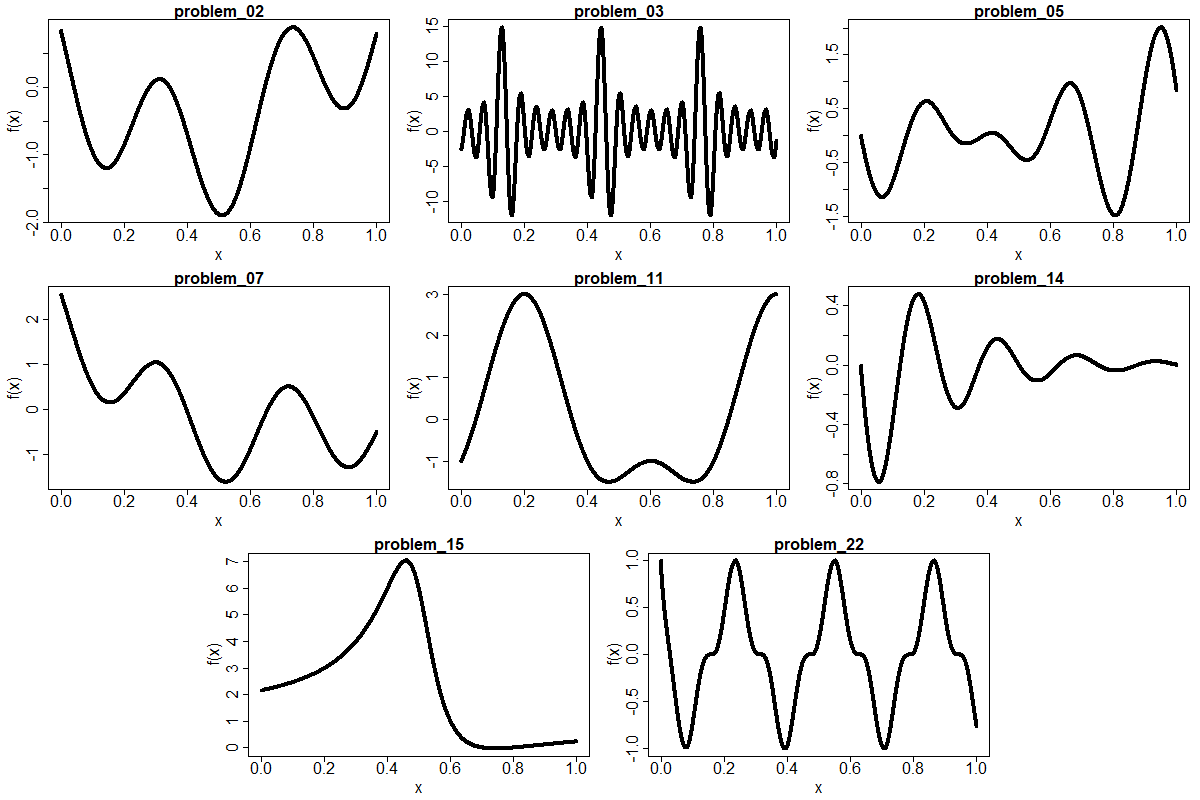}
    \caption{1-dimensional test problems. Original search spaces rescaled in $[0,1]$ }
    \label{fig:apx1}
\end{figure}

\newpage 

\begin{figure}
    \centering
    \includegraphics[width=\linewidth]{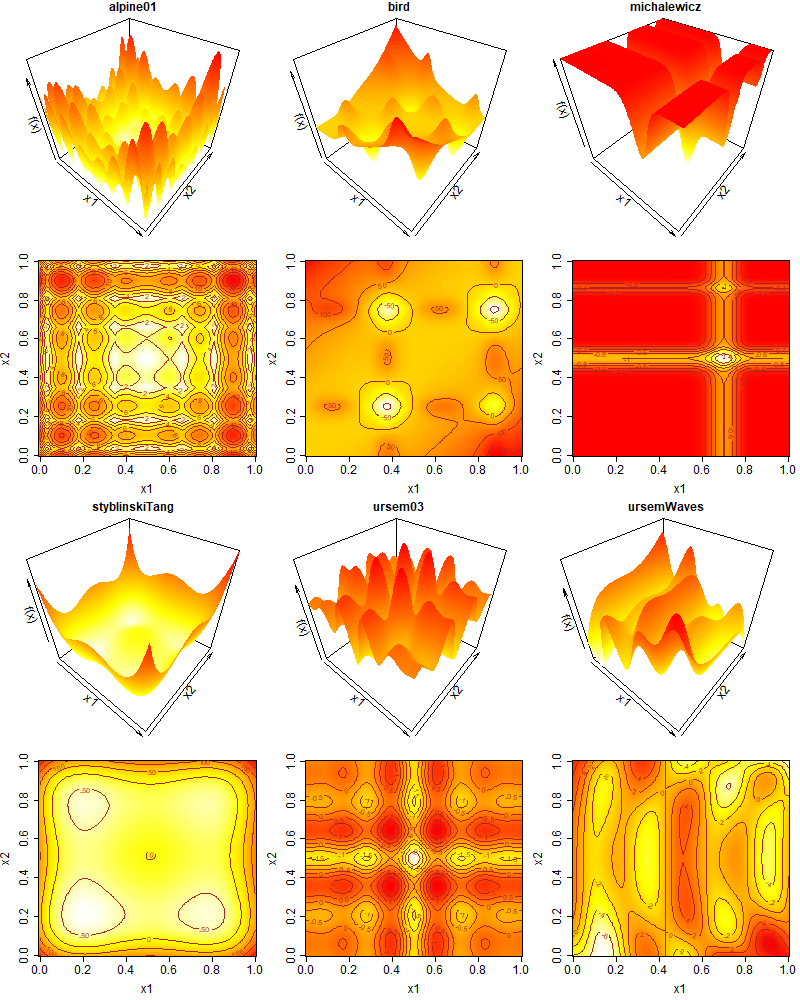}
    \caption{2-dimensional test problems (3D representation and countourplot for each test problem). Original search spaces rescaled in $[0,1]\times[0,1]$}
    \label{fig:apx2}
\end{figure}
\end{document}